\magnification=\magstep1
\font\rom=cmr10
\font\big=cmr12
\font\srom=cmr6
\font\goth=eufb10
\font\gras=cmbx12
\font\sit=cmti8

\def\square{\sqcap\hskip-0.228cm\sqcup}
\def\usquare{\sqcap\hskip-0.183cm\sqcup}
\def\nsquare{\sqcap\hskip-0.31cm\sqcup}
\def\msquare{\sqcap\hskip-0.23cm\sqcup}

\def\rsquare{\sqcap\hskip-0.255cm\sqcup}

\centerline{\gras On a symmetric space attached to polyzeta values.}

\vskip1cm

\centerline {\big Olivier Mathieu}

\vskip2cm
\centerline{ABSTRACT}
Quickly converging series are given to compute polyzeta numbers
$\zeta(r_1,\dots,r_k)$. The formulas involve an intricate combination
of (generalized) polylogarithms  at $1/2$. 
However, the combinatoric has a very simple 
geometric interpretation: it corresponds with the map $p\mapsto p^2$ on
a certain symmetric space $P$.

\footnote{}{

}

\vskip2cm
{\bf Introduction:} 

 Let $k\geq 1$. For a $k$-uple
$(r_1,$ $r_2$, $\dots,$ $r_k)$ of positive integers,
set  $\zeta(r_1,\dots,r_k)=
\sum_{0<n_1\dots<n_k}\, 1/n_1^{r_1}...n_k^{r_k}$.
We have $\zeta(r_1,\dots,r_k)<\infty$ if and only if $r_k\geq 2$.
By definition, the {\it polyzeta values} are the ${\bf Q}$-linear 
combinations of the finite numbers $\zeta(r_1,\dots,r_k)$.

Using the  definition of $\zeta(r_1,\dots,r_k)$,  the evaluation of a
polyzeta value  up to the $N^{\hbox{\rom th}}$ digit requires  to take into
account something like $O(10^N)$ terms. Therefore it is a very slow
computation.  A similar computational problem arises with the classical
series $\log 2=-\sum_{n>0} (-1)^n/n$ and
$\pi/4=\sum_{n\geq 0} (-1)^n 1/(2n+1)$, which converge very
slowly. 

However, we easily notice that:

\centerline{$\log 2=-\log (1-1/2)=\sum_{n>0} 2^{-n}/n$.}

\noindent A remarkable series for $\pi$ has been discovered by 
Bailey, Borwein and Plouffe {\bf [BBP]}:

\centerline{ 
$\pi=\sum_{n\geq 0} 1/2^{4n}[4/(8n+1) -2/(8n+4)-1/(8n+5)-1/(8n+6)]$}

\noindent Now to evaluate 
$\log 2$ or $\pi$ up to the $N^{\hbox{\rom th}}$ digit, one only 
needs the first $O(N)$-terms of the series and therefore 
$\log 2$ and $\pi$ can be computed very quickly. 
The goal of the paper is to 
provide similar identities  for all polyzeta values.

To do so, one needs to use the functions 
$L_{r_1,\dots, r_k}(z)= \sum_{0<n_1\dots<n_k}\, 1/n_1^{r_1}...n_k^{r_k}
z^{n_k}$,
 where $r_1,\dots, r_k$ are positive integers. By definition,
a ${\bf Q}$-linear  combinations of the functions
$L_{r_1,\dots, r_k}(z)$ is called 
a {\it polylogarithmic function}. The obvious identity
$\zeta(r_1,\dots,r_k)=L_{r_1,\dots, r_k}(1)$ does not help
to quickly evaluate polyzeta values.
However, the series defining polylogarithms  at $1/2$
converges very quickly:  to evaluate 
$L_{r_1,\dots, r_k}(1/2)$ up to the $N^{\hbox{\rom th}}$ digit, one only
needs to sum  $O(N^k)$-terms, and this can be done in polynomial time.
This remark suggests the 
following result:

\bigskip
MAIN STATEMENT: {\it Any polyzeta value is the value  of a
certain polylogarithmic function at $1/2$.}

\bigskip
In order to get a useful statement,  
the corresponding polylogarithmic function is described 
explicitly: see Theorem 7  for a precise
statement.  At first glance,  the combinatorics involved in Theorem 7 looks
intricated and therefore no details are given in the introduction. 
However, we can precisely formulate the main statement  in terms of
very simple geometric notions.  

Let $F(2)$ be the free group on two generators $\alpha$
and $\beta$ and let $s$ be the involution exchanging the generators. Let
$\Gamma={\bf Q}\otimes F_2$ be the Malcev completion of $\Gamma$ (see Section
4 for an alternative definition of $\Gamma$). The involution $s$ extends
to $\Gamma$ and there is a decomposition
$\Gamma=P.K$ where $K$ is the subgroup of fixed points of $s$ and
where $P=\{g\in\Gamma\vert s(g)=g^{-1}\}$. 
The group $\Gamma$ is proalgebraic over {\bf Q} and the symmetric space $P$
is a pro-algebraic variety over {\bf Q}. 

In section (4.8), all polyzeta values  are 
naturally indexed by rational functions on 
 $P$. Similarly, some 
polylogarithmic functions are naturally
indexed by rational functions on 
 $P$. So for $\phi\in{\bf Q}[P]$, denote
by $\zeta(\phi)$ and $L_{\phi}(z)$ the corresponding
polyzeta value and polylogarithmic function.

Now the square map $\msquare: P\rightarrow P, p\mapsto p^2$  induces an
algebra morphism $\msquare: {\bf Q}[P]\rightarrow {\bf Q}[P]$. The geometric
formulation of the main result is as follows:

\bigskip
MAIN THEOREM: {\it For any $\phi\in {\bf Q}[P]$,
$\zeta(\phi)=L_{\usquare\phi}(1/2)$. }

\bigskip
\noindent We also express polyzeta values as values of polylogarithmic
functions at $\rho^{\pm1}=\exp \pm i\pi/3$. The geometric interpretation
of this case is a bit more complicate because it involves an order $3$
automorphism of $\Gamma$, see section 4, Theorem 18.

\smallskip
{\it Acknowlegements:} A special thank to Wadim Zudilin. 
Section 5 has been suggested by him.

\smallskip
{\it Summary:}

1. Polylogarithms  and polyzeta values.

2. Polylogarithmic function at $1/2$ and at $\rho^{\pm1}$.

3. Explicit expressions for $\zeta(r)$.

4. Geometric interpretation of Theorem 7.

5.  Other expressions for zeta values. 

6. Conclusion.

\bigskip
{\bf 1. Polylogarithms and polyzeta values.}

This section is devoted to main definitions and conventions. 
The definitions of {\it polyzeta values}
and  {\it polylogarithmic functions} are not standard:
see the subsections (1.14) for more comments. 
Moreover in this section  we adopt some
conventions to renormalize infinite quantities like $\zeta(1)$ or 
$\int_0^z dt/t$.

\smallskip
(1.1) Shuffles: For  $N\geq 0$, denote by $S_N$  the symmetric group, i.e.
the set of all bijections $\sigma:\{1,\dots,N\}\rightarrow\{1,\dots,N\}$. 
Given $n$ and $m$ two non-negative integers, let $S_{n,m}$
be the set of all  $\phi\in S_{n+m}$ such that $\phi$ is
increasing on the subset $\{1,\dots,n\}$ and on the subset
$\{n+1,\dots,n+m\}$. The elements of $S_{n,m}$ are called 
{\it schuffles}.

\smallskip
(1.2) Shuffle product:   Let ${\cal W}$ be the set of words into the letters
$a$ and $b$. By convention, ${\cal W}$  contains the empty word $\emptyset$.
Set ${\cal H}={\bf Q}{\cal W}$, i.e. ${\cal H}$ is the ${\bf Q}$-vector
space with basis ${\cal W}$. For any two words $w=x_1\dots x_n$ and
$w'=x_{n+1}\dots x_{n+m}$, where each $x_i\in\{a,b\}$ is a letter, define
the product $w*w'\in{\cal H}$ by the formula:

$$w*w'=\sum_{\sigma\in S_{n,m}}\,x_{\sigma(1)}\dots x_{\sigma(n+m)}$$

\noindent By convention, we have $\emptyset *w=w*\emptyset=w$ for all word
$w$. The product $*$ is called the shuffle product. With respect to this
product, $\cal H$ is a commutative associative algebra, and $\emptyset$
is its unit.

\smallskip
(1.3) Subalgebras of ${\cal H}$: Let ${\cal W}^+$ be the set of words whose
first letter is not $b$. Equivalently, a word $w$ belongs to ${\cal W}^+$ if
$w=\emptyset$ or if $w$ starts with $a$. Similarly, let 
${\cal W}^{++}$ be the set of words whose first letter is
not $b$ and the last letter is not $a$. Set 
${\cal H}^+={\bf Q}{\cal W}^+$ and 
${\cal H}^{++}={\bf Q}{\cal W}^{++}$. It is easy to prove that
${\cal H}^+$ and  ${\cal H}^{++}$ are subalgebras of ${\cal H}$.

\smallskip
LEMMA 1: {\it There are isomorphisms of algebras:
${\cal H}={\cal H}^+[b]$ and 
${\cal H}={\cal H}^{++}[a,b]$.}

\smallskip
{\it Proof:} For each $n\geq 0$, let ${\cal W}_n$ be the set of words 
of the form $b^nw$ with $w\in {\cal W}^+$ , and set 
${\cal H}_n=\oplus_{0\leq k\leq n}\, {\bf Q}{\cal W}_i$.   We have 
$b*{\cal H}_n\subset {\cal H}_{n+1}$. Moreover we have
$b*w=(n+1) bw$  modulo ${\cal H}_n$ for any
$w\in {\cal W}_n$. It follows easily by induction that
${\cal H}_n=\oplus_{0\leq k\leq n}\,{\cal H}^{+}*b^k $,
i.e. ${\cal H}_n$ is the space of all polynomials in $b$ wih
coefficients in ${\cal H}^+$ and degree $\leq n$. Therefore
the first assertion follows.

   The proof of the second assertion is similar.

\smallskip
(1.4) The bijection $\lambda:{\cal W}^+\rightarrow\Lambda$: 

Let ${\bf N}$ be the set of positive integer.
For clarity, a word into the letters $1,2,\dots\in{\bf N}$ will
be called a {\it sequence of positive integers}.
Let $\Lambda$ the set of sequence  $(r_1\dots r_k)$ of positive integers. 
By convention, $\Lambda$  contains the empty sequence $\emptyset$.

Any word $w\in {\cal W}^+$ can be uniquely factorized as:
$w=ab^{t_1}ab^{t_2}\dots ab^{t_k}$, where $k$ is the number of
occurence of $a$ in $w$ and where the $t_i$ are non-negative integers.
Then, the map $w\in {\cal W}^+\mapsto (1+t_1,1+t_2,\dots,1+t_k)\in\Lambda$
defines a natural bijection $\lambda:{\cal W}^+\rightarrow\Lambda$.

\smallskip
(1.5) Polylogarithmic functions and polyzeta values:

Let $k\geq 1$ and let $r_1\dots r_k$ be a sequence of $k$ positive
integers. Consider the following  series in  the complex variable $z$:

$$L_{r_1,\dots,r_k}(z)=
\sum_{0<n_1<\dots <n_k}\,n_1^{-r_1}\dots n_k^{-r_k}\,z^{n_k}$$
In the infinite sum, the indices $n_1,\dots n_k$ are integers.
The functions
$L_{r_1,\dots,r_k}(z)$ are called  {\it polylogarithms}. 
Set $D=\{z\in{\bf C}\vert\,\vert z\vert<1\}$, 
$\overline{D}=\{z\in{\bf C}\vert\,\vert z\vert\leq 1\}$.
 The two points of interest for the
paper  are the following:

(i) if $r_k\geq 2$,  the series is
absolutely convergent on  ${\overline D}$ and therefore
$L_{r_1,\dots,r_k}(z)$ extends to a continous function on ${\overline D}$.

(ii) if $r_k=1$, the  series converges on $D$ and $L_{r_1,\dots,r_k}(z)$
extends to a continous function on ${\overline D}\setminus \{1\}$. 

\noindent Indeed
$L_{r_1,\dots,r_k}(z)$ extends to a
multivalued function,  see e.g. {\bf[C]} and  Proposition 2 below.
 For $r_k\geq 2$, set

$$\zeta(r_1,\dots,r_k)=
\sum_{0<n_1<\dots <n_k}\,n_1^{-r_1}\dots n_k^{-r_k}$$

In the paper, the numbers $\zeta(r_1,\dots,r_k)$ will be called
polyzeta values. Indeed the polyzeta value is  both the value at $z=1$
of the  polylogathm $L_{r_1,\dots,r_k}(z)$ and 
a  value of the polyzeta function
$\zeta(s_1,\dots,s_k)=
\sum_{0<n_1<\dots <n_k}\,n_1^{-s_1}\dots n_k^{-s_k}$.

\smallskip
(1.6) New notations:
Let $w\in{\cal W}^+$ and set $(r_1,\dots,r_k)=\lambda(w)$.
It is convenient to denote the function $L_{(r_1,\dots,r_k)}(z)$
by $L_w(z)$. Similarly set
$\zeta(w)=\zeta(r_1,\dots,r_k)$
if  $w\in{\cal W}^{++}$.

\smallskip
(1.7) The one-forms $\omega_a$ and $\omega_b$:
Define the following one-forms on {\bf C}:

$$\omega_a(z)={\hbox
{\rom d}z\over1-z}\hskip0.4cm\hbox{\rom and}\hskip0.4cm
\omega_b(z)={\hbox{\rom d}z\over z}$$

\noindent For an element $c=xa+yb\in {\bf
Q}a\oplus{\bf Q}b$, set
$\omega_c(z)=x\omega_a(z)+y\omega_b(z)$.
Given a  smooth path
$\gamma:[0,1]\rightarrow{\bf C}, t\mapsto\gamma(t)$, recall that
$\gamma^*\omega_a(t)= {\gamma'(t)
\over1-\gamma(t)}\,\,\hbox{\rom d}t$, $\gamma^*\omega_b(t)= {\gamma'(t)
\over\gamma(t)}\,\,\hbox{\rom d}t$.  and
$\gamma^*\omega_c(t)=x\gamma^*\omega_a(t)+y\gamma^*\omega_b(t)$.

\smallskip
(1.8) Kontsevitch formula:
For a positive integer $n$, set
$\Delta_n=\{(x_1, x_2\dots, x_n)\in{\bf R}^n\vert
0\leq x_1\leq x_2\dots\leq x_n\leq1\}$.
Let $w=c_1\dots c_n\in {\cal W}^{+}$ be a word, where each $c_i$ is
a letter.  The following formula is due to Kontsevitch (see {\bf [Z]}).

\smallskip
PROPOSITION 2: {\it Let $w=c_1\dots c_n\in {\cal W}^{+}$ be a word, let 
$z\in \overline{D}$,  and let
$\gamma:[0,1]\rightarrow
\overline{D}$ be a path with
$\gamma(0)=0$ and $\gamma (1)=z$. 

Assume that $w\in {\cal W}^{++}$ or that $\gamma$ does not meet $1$. 
Then we have:

$$L_{w}(z)=\int_{\Delta_n}\,\gamma^*\omega_{c_1}(x_1)
\gamma^*\omega_{c_2}(x_2)\dots\gamma^*\omega_{c_n}(x_n)$$.}

\smallskip
In {\bf [Z]}, Kontsevitch formula is stated for the straight path
$t\mapsto zt$, but it is easy to see
that the integral is homotopy invariant as long $\gamma$ stay in 
$\overline D$ (and $\gamma$ stay  $\overline D\setminus\{1\}$ if
$w\notin{\cal W}^{++}$).

\smallskip
(1.9) Products: Let $n,m$ be non negative
integers and let $c_1,c_2,\dots,c_{n+m}\in \{a,b\}$ be letters
with $c_1=c_{n+1}=a$.
Set  $u=c_1\dots c_n$ and  $v=c_{n+1}\dots c_{n+m}$.  For 
$\sigma\in S_{n,m}$, set 
$w_{\sigma}=c_{\sigma(1)}\dots c_{\sigma(n+m)}$.

\smallskip
COROLLARY 3:
{\it 
For $u,\,v\in {\cal W}^{+}$, we have $L_u(z)L_v(z)=
\sum_{\sigma\in S_{n,m}}\,\,L_{w_{\sigma}}(z)$ for all $z\in D$.
Moreover for $u,\,v\in {\cal W}^{++}$, we have
$\zeta(u)\zeta(v)=\sum_{\sigma\in S_{n,m}}\,\zeta(w_{\sigma})$.}

\smallskip
{\it Proof:} 
Set 

\centerline{$\omega'=
\omega_{c_1}(zx_1)\omega_{c_2}(zx_2)\dots\omega_{c_n}(zx_n)$, }

\centerline{$\omega"= \omega_{c_{n+1}}(zx_{n+1})
\omega_{c_{n+2}}(zx_{n+2})\dots\omega_{c_{n+m}}(zx_{n+m})$,}

\centerline{$\Delta_m= \{(x_{n+1}, x_{n+2}\dots, x_{n+m})\in{\bf R}^n
\vert 0\leq x_1\leq x_2\dots\leq x_n\leq1\}$,}

\noindent and for  $\sigma\in S_{n,m}$, set 

\centerline{$\Delta_{\sigma}=\{(x_1,\dots,x_{n+m})\in{\bf R}^{n+m}\vert
x_{\sigma(1)}\leq x_{\sigma(2)}\leq \dots\leq x_{\sigma(n+m)}\}$.}

\noindent Since  $\Delta_{n}\times\Delta_m=
\cup_{\sigma\in S_{n,m}}\,\Delta_{\sigma}$, we get
$\int_{\Delta_n}\omega'\,\int_{\Delta_m}\omega"=
\int_{\Delta_n\times\Delta_m}\omega'\wedge\omega"=
\sum_{\sigma}\,\int_{\Delta_{\sigma}}\,\omega'\wedge\omega"$.
By Proposition 2, this identity is equivalent to
$L_u(z)L_v(z)=
\sum_{\sigma\in S_{n,m}}\,\,L_{w_{\sigma}}(z)$.
At $z=1$, one gets the second identity
$\zeta(u)\zeta(v)=\sum_{\sigma\in S_{n,m}}\,\zeta(w_{\sigma})$. Q.E.D.

\smallskip
(1.10) Final definitions and notations for polylogarithmic functions:
Up to now, the polylogarithms $L_w(z)$ are 
defined for $w\in {\cal W}^{+}$. In order to extend the 
definition to all $w\in {\cal W}$,  a renormalization 
procedure is used.

Set $\Omega=D\setminus ]-1,0]$
and let $Hol(\Omega)$ be the algebra of holomorphic functions on
$\Omega$. Since $\omega$ is simply connected, let denote by $\log z$ the
holomorphic function on $\Omega$ whose  restriction to $]0,1[$ is the usual
logarithmic function.

\smallskip
LEMMA 4:
{\it There is  a unique algebra morphism
$\Phi:{\cal H}\rightarrow Hol(\Omega)$ such that 
$\Phi(w)=L_w(z)$ for $w\in{\cal H}^+$ and
$\Phi(b)=\log z$. }

\smallskip
{\it Proof:} This follows from Lemmma 1 and corollary 3. Q.E.D.

\smallskip
For any $h\in {\cal H}$, set $L_h(z)=\Phi(h)$.
When $h$ is a word $w$ in ${\cal W}^+$, 
this new notation agrees with 
the previous one.  Corollary 3 can be restated as:
$L_u(z)L_v(z)= L_{u*v}(z)$.
By definition the {\it polylogarithmic functions} are the functions
$L_h(z)$ whith $h\in{\cal H}$. 

This definition is a
slighty different from the introduction. 
However, we will only use polylogarithmic functions
$L_h(z)$ with $h\in{\cal H}^{+}$, which are the polylogaritmic 
functions defined in introduction.

\smallskip
(1.11) Final definitions and notations for polyzeta values:
Up to now, the polyzeta values $\zeta(w)$ are 
defined for $w\in {\cal W}^{++}$. In order to extend the 
definition to all $w\in {\cal W}$,  we will use a renormalization 
procedure as follows.

\smallskip 
LEMMA 5:
{\it There are three algebra morphisms
$\psi,\,\psi^+,\,\psi^-: {\cal H}\rightarrow {\bf C}$ 
uniquely defined by the following requirements:

$\psi(w)=\psi^+(w)=\psi^-(w)=\zeta(w)$ if $w\in{\cal H}^{++}$

$\psi(a)=0$, $\psi^+(a)=i\pi$, $\psi^-(a)=-i\pi$

$\psi(b)=\psi^+(b)=\psi^-(b)=0$.}

\smallskip
{\it Proof:} This follows from Lemmma 1 and corollary 3.

\smallskip
Similarly, this allows to define $\zeta(h)=\psi(h)$,
$\zeta^{\pm}(h)=\psi^{\pm}(h)$ for any
$h\in{\cal H}$. By definition the {\it polyzeta values} are the 
numbers $\zeta(h)$ whith $h\in{\cal H}$. 
Corollary 3 can be restated as: $\zeta(u)\zeta(v)=\zeta(u*v)$
and $\zeta^{\pm}(u)\zeta^{\pm}(v)=\zeta^{\pm}(u*v)$ for any
$u,\,v\in{\cal H}$.

Set  ${\cal Z}=\zeta({\cal H})$ and ${\cal Z}^{\pm}=\zeta({\cal H}^{\pm})$.
By definition, ${\cal Z}$ and ${\cal Z}^{\pm}$ are subrings of ${\bf C}$,
and ${\cal Z}$ is the space of all polyzeta values. It is easy to compare
the three algebras ${\cal Z}$ and ${\cal Z}^+$ and ${\cal Z}^-$.

\smallskip
LEMMA 6: {\it 

(i) As a ${\bf Q}$ vector space, ${\cal Z}$ is generated by
all $\zeta(w)$ with $w\in {\cal W}^{++}$.

(ii) We have ${\cal Z}\subset {\bf R}$.

(iii) ${\cal Z}^{\pm}={\cal Z}\oplus i\pi {\cal Z}$.}

\smallskip
{\it Proof:} The assertions (i) and (ii) follow from Lemma 1.
Moreover it follows that ${\cal Z}^{\pm}$ is the {\bf Q}-algebra generated
by ${\cal Z}$ and $\zeta^{\pm}(b)=\pm i\pi$. However
 $(\zeta^{\pm}(b))^2=-\pi^2=-6\zeta(2)$, therefore 
$(\zeta^{\pm}(b))^2$ belongs to ${\cal Z}$ and assertion (6.3) follows.

\smallskip
(1.12) Hopf algebra structure:
Define the linear maps
$\eta: {\cal H}\rightarrow{\bf Q}$, 
$\iota:{\cal H}\rightarrow{\cal H}$ and
$\Delta:{\cal H}\rightarrow{\cal H}\otimes{\cal H}$
as follows.
For any word $w=c_1\dots c_n\in{\cal W}$, set

$$\eta(w)=1\, \,\hbox{\rom if}\, w=\emptyset\,\, \hbox{\rom and}\, \eta(w)=0
\,\,\hbox{\rom otherwise}$$

$$\iota(w)=(-1)^n\,c_n c_{n-1}\dots c_n\ $$
 
$$\Delta(w)=\sum_{0\leq i\leq n}\,c_1\dots c_i\otimes c_{i+1}\dots c_n$$

\noindent
The map $\eta$, $\iota$ and $\Delta$ are algebra morphisms. 
Indeed ${\cal H}$ is a Hopf algebra  with co-unit $\eta$,
inverse map $\iota$ and  coproduct $\Delta$.

\smallskip
(1.13) Concatenation product:
For two words $w=c_1\dots c_n$ and $w'=c_{n+1}\dots c_{n+m}$,
their concatenation is the word $w w'=c_1\dots c_n c_{n+1}\dots
c_{n+m}$. This induces another structure of algebra on ${\cal H}$,
for which the product of two elements $h$, $h'$ is simply denoted by 
$h\,h'$.

\smallskip
(1.14) Remarks on references and on the terminology:
  
In the classical litterature, only the
functions $L_k(z)=\sum_{n>0} z^n/n^k$ are called polylogarithms, see {\bf
[L]} {\bf [Oe]}.  We did not find a standard name for the
$L_{(r_1,\dots,r_k)}(z)$. They are defined in the Bourbaki's talk {\bf [C]},
where the title  suggests to call them again polylogarithms.

It seems that some polyzeta values, like $\zeta(1,3)$, were already
known by  Euler, see {\bf [C]}. The general definition of 
$\zeta(r_1,\dots,r_k)$ appears explicitely around 1990 in {\bf [H]} and {\bf
[Z]}. These numbers are also called multiple zeta values in {\bf [Z]},
multiple harmonic sums in {\bf [H]}, multizeta numbers in {\bf [E]},
Euler-Zagier numbers in {\bf [BB]} and polyzetas numbers in {\bf[C]}.

The fact that polyzeta values are naturally indexed by words
has been observed by many authors,  see {\bf [H]},
{\bf [H-P]} and {\bf [C]}.  Lemma 1 and corollary 3 are well-known. Proofs
are given for the  convenience of the reader.

\vskip1cm
{\bf 2. Polylogaritmic functions at 1/2 and at $\rho^{\pm1}$:}

Define two linear maps $\sigma, \,\tau:{\cal H}\rightarrow {\cal H}$
as follows. First set $\sigma(a)=b$, $\sigma(b)=a$,
$\tau(a)=a+b$ and $\tau(b)=-a$.  For any word
$w=c_1\dots c_n\in {\cal W}$, set
$\sigma(w)=\sigma(c_n)\dots\sigma(c_1)$ and 
$\tau(w)=\tau(c_n)\dots\tau(c_1)$. It is easy to see that $\sigma$ and
$\tau$ are algebra morphisms relative to the shuffle product $*$
(they are anti-morphism relative to the concatenation product).

Define now the two operators 
${\msquare}, \nabla:{\cal H}\rightarrow {\cal H}$ as the
following composite maps:

$$\msquare:{\cal H}\buildrel\Delta\over\rightarrow
{\cal H}\otimes{\cal H}\buildrel id\otimes\sigma\over\longrightarrow
{\cal H}\otimes{\cal H}\buildrel *\over\rightarrow
{\cal H}$$

$$\nabla: {\cal H}\buildrel\Delta\over\rightarrow
{\cal H}\otimes{\cal H}\buildrel id\otimes\tau\over\longrightarrow
{\cal H}\otimes{\cal H}\buildrel *\over\rightarrow
{\cal H}$$

Set $\rho=e^{i\pi/3}$.

\bigskip
THEOREM 7:
{\it For any  $h\in{\cal H}$, we have:

$$\zeta(h)=L_{\usquare(h)}(1/2)$$

$$\zeta^{\pm}(h)=L_{\nabla(h)}(\rho^{\pm1})$$}

\bigskip
{\it Proof:} 
Note that $\rsquare(a)=\rsquare(b)= a+b$,
$L_a(z)=\log 1/(1-z)$ and $L_b(z)=\log z$, therefore
$L_a(1/2)+ L_b(1/2)=0=\zeta(a)=\zeta(b)$. Similarly,
$\nabla(a)=2a+b$ and $\nabla(b)=0$, and we have
$L_\nabla(a)(\rho^{\pm 1})=
-2\log (1-\rho^{\pm 1})+ \log(\rho^{\pm 1})=\pm i\pi=\zeta^{\pm}(a)$
and $L_\nabla(b)(\rho^{\pm 1})=0=\zeta^{\pm}(b)$.

Since $\msquare$ and $\nabla$ are algebra morphisms, and since 
the algebra ${\cal H}$ is generated by $a$, $b$ and
${\cal W}^{++}$, it is enough to show 
the formulas when $h$ is a non empty word $w$ in ${\cal W}^{++}$.
So let  $w\in{\cal W}^{++}$ be a word of length $n\geq 2$.
Set $w=c_1\dots c_n$ where
$c_i\in\{a,\,b\}$ are letters with $c_1=a$ and $c_n=b$. 

Let $\gamma:[0,1]\rightarrow \overline{D},\,t\mapsto t$ be the straight path
from $0$ to $1$. Choose two smooth paths $\,\gamma_{\pm}:
[0,1]\rightarrow \overline{D}$ with the following properties:
$\,\gamma_{\pm}(0)=0$, $\,\gamma_{\pm}(1)=1$, 
$\,\gamma_{\pm}(1/2)=\rho^{\pm1}=1/2\pm i\sqrt{3}/2$ and 
$\hbox{\rom Re}\, \gamma_{\pm}(t)\geq 1/2$ for all $t\in[1/2,1]$.
Set $\Delta_n=\{(x_1, x_2\dots, x_n)\in{\bf R}^n\vert
0\leq x_1\leq x_2\dots\leq x_n\leq1\}$. By Proposition 2, we have:

$$\zeta(w)=\int_{\Delta_n}\,\eta^*\omega_{c_1}(x_1)\,
\eta^*\omega_{c_2}(x_2)\dots\eta^*\omega_{c_n}(x_n)$$

\noindent where $\eta$ is the path $\gamma$, or $\gamma^+$ or
$\gamma^-$.

 For $0\leq i\leq n$,
set  $\Delta'_i=\{(x_1, x_2\dots, x_i)\in{\bf R}^n\vert
0\leq x_1\leq x_2\dots\leq x_i\leq 1/2\}$ and 
$\Delta"_i=\{(x_{i+1}, x_{i+2}\dots, x_n)\in{\bf R}^n\vert
1/2 \leq x_{i+1}\leq x_{i+2}\dots\leq x_n\leq1\}$. 
From the decomposition:
$\Delta=\cup_{0\leq i\leq n}\, \Delta'_i\times \Delta"_i$,
it follows that $\zeta(w)=\sum_{0\leq i\leq n} L'_i\, L"_i
=\sum_{0\leq i\leq n} L_i'^{\pm}\, L"_i^{\pm}$,
where the numbers $L'_i,\, L"_i$,
$L_i'^{\pm},\, L_i"^{\pm}$ are the following integrals:

$$L'_i=  \int_{\Delta'_i}\,\gamma^*\omega_{c_1}(x_1)\,
\gamma^*\omega_{c_2}(x_2)\dots\gamma^*\omega_{c_i}(x_1)$$

$$L"_i=  \int_{\Delta"_i}\,\gamma^*\omega_{c_{i+1}}(x_{i+1})\,
\gamma^*\omega_{c_{i+2}}(x_{i+2})\dots\gamma^*\omega_{c_n}(x_n)$$

$$L_i'^{\pm}=  \int_{\Delta'_i}\,\gamma_{\pm}^*\omega_{c_1}(x_1)\,
\gamma_{\pm}^*\omega_{c_2}(x_2)\dots\gamma_{\pm}^*\omega_{c_i}(x_1)$$

$$L"_i^{\pm}=  \int_{\Delta"_i}\,\gamma_{\pm}^*\omega_{c_{i+1}}(x_{i+1})\,
\gamma_{\pm}^*\omega_{c_{i+2}}(x_{i+2})\dots\gamma_{\pm}^*\omega_{c_n}(x_n)$$

\noindent
Using Kontsevitch formula, we get $L'_i= L_{w'_i}(1/2)$
and $L_i'^{\pm}=L_{w'_i}(\rho^{\pm 1})$, where 
$w'_i=c_1\dots c_i$. 
To evaluate $L"_i$, one needs to introduce some new notations.
Define by $S,\,T:{\bf C}\rightarrow{\bf C}$ the rational maps:
$S(z)=1-z$ and $T(z)=1-1/z$. Define the new paths 
$\delta,\,\delta_{\pm}:[0,1/2]\rightarrow{\bf C}$ by
$\delta(t)=1-\gamma(1-t)=S\circ\gamma(1-t)$ and 
$\delta_{\pm}(t)=1-1/\gamma_{\pm}(1-t)=T\circ\gamma(1-t)$.  
Clearly, $\delta$ is the straight path from $0$ to $1/2$.
Since $T(\rho^{\pm1})=\rho^{\pm1}$ and $T(1)=0$, 
$\delta_{\pm}$ is a path from $0$ to
$\rho^{\pm1}$. Since 
$\hbox{\rom Re}\, \gamma_{\pm}(t)\geq 1/2$ for all $t\in[1/2,1]$,
it follows that  $\delta_{\pm}$ lies in $\overline {D}\setminus\{1\}$.

With the
convention of (1.7), we get:

$$\gamma^*\omega_c(t)=\delta^*\omega_{\sigma(c)}(1-t)\,\,\hbox{\rom and}
\,\,\gamma_{\pm}^*\omega_c(t)=\delta_{\pm}^*\omega_{\tau(c)}(1-t)$$

\noindent for any $c\in{\bf Q}a\oplus{\bf Q}b$.
Using the new variables $y_j= 1-x_j$,
we thus get:

$$L"_i=  \int_{\overline{\Delta"_i}}
\,\delta^*\omega_{\sigma(c_{n})}(y_{n})
\delta^*\omega_{\sigma(c_{n-1})}(y_{n-1})\dots
\delta^*\omega_{\sigma(c_{i+1})}(y_{i+1})$$

$$L"_i^{\pm}=  \int_{\overline{\Delta"_i}}
\,\delta_{\pm}^*\omega_{\tau(c_{n})}(y_{n})
\delta_{\pm}^*\omega_{\tau(c_{n-1})}(y_{n-1})\dots
\delta_{\pm}^*\omega_{\tau(c_{i+1})}(y_{i+1})$$

\noindent where $\overline{\Delta"_i}=
\{(y_{n}, y_{n-1},\dots, y_{i+1})\in{\bf R}^n\vert
0 \leq y_n\leq y_{n+1}\dots\leq y_{i+1}\leq 1/2\}$. 
It follows from Proposition 2 that $L"_i= L_{\sigma(w"_i)}(1/2)$
and $L"_i^{\pm}=L_{\tau(w"_i)}(\rho^{\pm1})$
where $w"_i$ is the word $c_{i+1}\dots c_n$. Therefore we get

$$\zeta(w)=\sum_{0\leq i\leq n}
L_{w'_i}(1/2)  L_{\sigma(w"_i)}(1/2),\,\hbox{\rom and}$$

$$\zeta(w)=\sum_{0\leq i\leq n}
L_{w'_i}(\rho^{\pm1})  L_{\tau(w"_i)}(\rho^{\pm1}).$$

\noindent Since $\Delta(w)=\sum _{0\leq i\leq n}\, w'_i \otimes w"_i$,
it is clear that 
$\nsquare (w)=\sum_{0\leq i\leq n}\, w'_i * \sigma(w"_i)$
and $\nabla(w)=\sum_{0\leq i\leq n}\, w'_i * \tau(w"_i)$, 
and therefore the formula follows from Corollary 3. Q.E.D.

\vskip1cm
{\bf 3. Explicit expressions for $\zeta(r)$.}

Theorem 7 provides a combinatorial way to express any polyzeta value
as  a polylogarithmic functions at $1/2$
or  at $\rho$ or at $\overline{\rho}$. In this section,
Theorem 10 and Corollary 12 provided closed formulas for zeta
values $\zeta(r)$, where $r\geq 2$ is a given integer. 
The formulas are derived from Theorem 7. 
However, for general polyzeta values 
$\zeta(r_1,\dots,r_k)$ the combinatorics seem too intricate to find a
simple combinatorial formula.

The concatenation product hh', which is not commutative, should not be
confused with the commutative schuffle product $h*h'$. The following
conventions will be used.
First, for $h\in{\cal H}$ and $n\geq 1$, the notation $h^n$
will be the $n^{\hbox{\sit th}}$ power of $h$ with respect to the
concatenation product. Moreover, the concatenation product takes precedence
of the schufle product. For example, the  expression $hh'*h"$ should be
understood as
$(hh')*h"$.

\smallskip
{\it LEMMA 8:} {\it We have:
 
$$\nsquare(ab^{r-1})=
2 a^{r-1}(a+b) +\sum_{1\leq j\leq r-2} a^j(a+b)^{r-j}$$

$$\nabla(ab^{r-1})=
(-1)^{r+1} 3 a^r +\sum_{1\leq j\leq r-2} (-1)^{j+1} a^j(b-a)^{r-j}$$}

\smallskip
{\it Proof:} 
We have:

$\nsquare(ab^{r-1}) =\sigma(ab^{r-1}) + \sum\limits_{0\leq j\leq
r-1}\,ab^{j}*\sigma(b^{r-1-j})$

\hskip1.6cm $=a^{r-1}b + \sum\limits_{u+v= r-1}\,ab^{u}*a^{v}$.

 Similarly, we have:

$\nabla(ab^{r-1}) =\tau(ab^{r-1}) + \sum\limits_{0\leq j\leq
r-1}\,ab^{j}*\tau(b^{r-1-j})$

\hskip1.6cm $=(-1)^{r+1}a^{r-1}(a+b) + \sum\limits_{u+v= r-1}\,(-1)^v
ab^{u}*a^{v}$.

Use now the formula, 
$cw*a^v=\sum\limits_{i+j=v}\,a^jc(w*a^i)$, which holds for any
word $w$ and any letter $c$. Thus we have
$ab^{u}*a^{v}=\sum\limits_{i+j=v}\,a^{j+1}(b^u*a^i)$ and  we get:

$\nsquare(ab^{r-1}) = a^{r-1}b + \sum\limits_{j+u+v=
r-1}\,a^{j+1}(b^{u}*a^{v})$, and

$\nabla(ab^{r-1}) = (-1)^{r+1}a^{r-1}(a+b) + 
\sum\limits_{j+u+v= r-1}\,(-1)^{j+v}\,a^{j+1}(b^{u}*a^{v})$.

Using now the formulas:

$(a+b)^N=\sum\limits_{u+v=N}\, a^u*b^v$ and 
$(b-a)^N=\sum\limits_{u+v=N}\, (-1)^u\,a^u*b^v$

we get:

$\nsquare(ab^{r-1}) = a^{r-1}b + \sum_{0\leq j\leq
r-1}\,a^{j+1}(a+b)^{r-1-j}$

\hskip1.6cm $=a^{r-1}b + \sum_{1\leq j\leq
r}\,a^j(a+b)^{r-j}$

\hskip1.6cm $=a^{r-1}b + a^r +\sum_{1\leq j\leq
r-1}\,a^{j}(a+b)^{r-j}$

\hskip1.6cm $=
 a^{r-1}(a+b) +\sum_{1\leq j\leq r-1} a^j(a+b)^{r-j}$

\hskip1.6cm $=
 2 a^{r-1}(a+b) +\sum_{1\leq j\leq r-2} a^j(a+b)^{r-j}$. 

We also get:

$\nabla(ab^{r-1}) = (-1)^{r+1}a^{r-1}(a+b) + 
+ \sum\limits_{0\leq j\leq
r-1}\,(-1)^j a^{j+1}(b-a)^{r-1-j}$

\hskip1.6cm $=(-1)^{r+1}a^{r-1}(a+b) +  
 \sum\limits_{1\leq j\leq r}\,(-1)^{j+1}a^j(b-a)^{r-j}$

\hskip1.6cm $=(-1)^{r+1}a^{r-1}(a+b) + (-1)^{r+1}a^r +
(-1)^ra^{r-1}(b-a)$

\hskip5cm$ +
\sum\limits_{1\leq
j\leq r-2}\,(-1)^{j+1}a^j(b-a)^{r-j}$

\hskip1.6cm 
$=(-1)^{r+1} 3 a^r +\sum_{1\leq j\leq r-2} (-1)^{j+1} a^j(b-a)^{r-j}$.
Q.E.D.

\smallskip
Let ${\cal W}_r$ be the set of words of lenght $r$.
Any $w\in {\cal W}_r$ can be written as 
$w=a^ju$, where $u$ does not start with $a$.
Let $k$ be the number of occurence of $a$ in $w$.
Set $j(w)=j$, $k(w)=k$ and define the numbers
$c(w)$ and $c^{\pm}(w)$ as follows.

(i) If $w=a^r$, set 
$c(w)=r$ and $c^{\pm}(w)=(-1)^{r+1}(r+1)$.

(ii) If $w=a^{r-1}b$, set
$c(w)=r$ and $c^{\pm}(w)=(-1)^{r}(r-2)$.

(iii) Otherwise,  set 
$c(w)=j$ and $c^{\pm}(w)=(-1)^{k+1}j$.

\smallskip
{\it LEMMA 9:} {\it We have:
 
$$\nsquare(ab^{r-1})=
\sum_{w\in{\cal W}_r} c(w) w$$

$$\nabla(ab^{r-1})=
\sum_{w\in{\cal W}_r} c^{\pm}(w) w$$}

\smallskip
{\it Proof:}
The first identity of Lemma 8 can be written as:

$$\nsquare(ab^{r-1})=
a^{r-1}b +\sum_{1\leq i\leq r} a^i(a+b)^{r-i}$$.

Since $(a+b)^{r-i}=\sum_{u\in{\cal W}_{r-i}}\, u$, we thus get

$$\nsquare(ab^{r-1})=
a^{r-1}b +\sum_{1\leq i\leq r}\sum_{u\in{\cal W}_{r-i}}\, a^iu$$

The word $w=a^{j(w)}v$ belongs to $a^i{\cal W}_{r-i}$ for
all $i\leq j(w)$. Therefore

$\sum_{1\leq i\leq r} a^i(a+b)^{r-i}=
\sum_{w\in{\cal W}_r}\, j(w) w$. Since 
$c(a^{r-1}b)=j(a^{r-1}b)+1$ and $c(w)=j(w)$ otherwise,
the formula $\nsquare(ab^{r-1})=
\sum_{w\in{\cal W}_r} c(w) w$ is now proved.

Set
$\overline a=-a$ and
$\overline b=b$. For a word
$w=c_1\dots c_n$, set $\overline w=\overline {c_1}\dots\overline {c_r}$
and for a general element $h=\sum c_w w$ in ${\cal H}$ set
$\overline h=\sum c_w \overline{w}$. Since the involution
$h\mapsto\overline{h}$ is a morphism relative to the concatenation
product, it follows from Lemma 8 that:

$-\overline{\nabla(ab^{r-1})}=
 3 a^r +\sum_{1\leq i\leq r-2}  a^i(b+a)^{r-i}$

\hskip1.8cm$=
  a^r-a^{r-1}b +\sum_{1\leq i\leq r}  a^i(b+a)^{r-i}$

It follows from the previous proof that
$\sum_{1\leq i\leq r} a^i(a+b)^{r-i}=
\sum_{w\in{\cal W}_r}\, j(w) w$. Therefore, one gets:

$-\overline{\nabla(ab^{r-1})}=
a^r-a^{r-1}b +\sum_{w\in{\cal W}_r}\, j(w) w$.

Note that $\overline w=(-1)^{k(w)}\,w$ for all words $w$. Thus:

$\nabla(ab^{r-1})= (-1)^{r+1} a^r - (-1)^ra^{r-1}b 
+\sum_{w\in{\cal W}_r}\, (-1)^{(1+k(w))}j(w) w$

Since 
$c^{\pm}(a^r)=(-1)^{r+1}(j(a^r+1)$,
$c^{\pm}(a^{r-1}b)=(-1)^r(j(a^{r-1}b)-1)$ and 
$c^{\pm}(w)=(-1)^{1+k(w)}\,j(w)$ otherwise,
the formula $\nabla(ab^{r-1})=
\sum_{w\in{\cal W}_r} c^{\pm}(w) w$ is now proved. Q.E.D.

\bigskip
Let $r\geq 2$ be an integer. Let $\Lambda_r$ be the set of all
${\bf m}=(m_1,\dots m_k)\in\Lambda$  with $m_1+\dots+m_k=r$. For ${\bf
m}=(m_1,\dots m_k)\in\Lambda_r$, set $k({\bf m})=k$ and define the integers 
$j({\bf m})$, $b({\bf m})$ and $b^{\pm}({\bf m})$ as follows.

(i) If $m_1=m_2=\dots=m_r=1$, set $j({\bf m})=r$,
$b({\bf m})=r$ and $b^{\pm}({\bf m})=(-1)^{r+1}(r+1)$.

\noindent Otherwise, let $j({\bf m})=j$ be the index such that
$m_1=m_2=\dots=m_{j-1}=1$ and $m_{j}\geq 2$.

(ii) If $m_1=m_2=\dots=m_{r-2}=1$ and
$m_{r+1}=2$, set
$b({\bf m})=r$ and $b^{\pm}({\bf m})=(-1)^{r}(r-2)$. 

(iii) Otherwise,  set 
$b({\bf m})=j$ and $b^{\pm}({\bf m})=(-1)^{k+1}j$,
where $j=j({\bf m})$ and $k=k({\bf m})$.

\bigskip
THEOREM 10: {\it For $r\geq 2$, we have:

$$\zeta(r)=\sum_{{\bf m}\in \Lambda_r}\,\, b({\bf m}) L_{\bf m}(1/2)$$

$$\zeta(r)=\sum_{{\bf m}\in \Lambda_r}\,\,
 b^{\pm}({\bf m}) L_{\bf m}(\rho^{\pm1})$$}

\bigskip

{\it Proof:} It is clear that  $c(w)$ and
$c^{\pm}(w)$ vanish if $w\notin{\cal W}^+$. 
Therefore it follows from Theorem 7 and Lemma 9 that

$$\zeta(r)=\sum_{w\in{\cal W}^+_r} c(w) L_w(1/2)$$

$$\zeta(r)=
\sum_{w\in{\cal W}^+_r} c^{\pm}(w) L_w(\rho^{\pm1})$$

\noindent where ${\cal W}^+_r={\cal W}^+\cap {\cal W}_r$.
Note that the map $\lambda$ of section 1.4 provides a
bijection $\lambda: {\cal W}^+_r\rightarrow \Lambda_r$.
It is easy to check that 
$j(w)=j(\lambda(w))$ and $k(w)=k(\lambda(w))$ for all
$w\in{\cal W}^+_r$, and therefore

$$c(w)=b(\lambda(w))\hskip.5cm\hbox{\rom and}\hskip.5cm
c^{\pm}(w)=b^{\pm}(\lambda(w))$$

\noindent for all $w\in{\cal W}^+_r$. Therefore Theorem 10 is proved.
Q.E.D.

\smallskip
For $1\leq i\leq r-1$, set 

$$C_i=\{{\bf n}=(n_1\dots n_k)\in{\bf Z}^r\vert 
0<n_1<...<n_i\leq n_{i+1}\leq\dots \leq n_k\}.$$

\noindent Also, set $C_r=C_{r-1}$.

\smallskip
LEMMA 11:
{\it We have:

$$\zeta(r)= \sum_{1\leq i\leq r}\,
\sum_{{\bf n}\in C_i} \,\,{2^{-n_r}\over n_1n_2\dots n_r }$$.}

\smallskip
{\it Proof:} Set $c=a+b$. For any word 
$w=d_1\dots d_r$ into the letters $a$, $b$ and $c$, let
$C_w$ be the set of all ${\bf n}=(n_1\dots n_k)\in{\bf Z}^r$ satisfying
the following property:

$$0{\cal R}_1 n_1 {\cal R}_2 n_2\dots {\cal R}_r n_r$$

\noindent where ${\cal R}_i$ stands for the symbol $<$ if $d_i=a$,
${\cal R}_i$ stands for  the symbol $=$ if $d_i=b$ and 
${\cal R}_i$ stands for  the symbol $\leq$ if $d_i=c$.
So if $w$ is a word  into the letters $a$, $b$ and $c$,
we get $L_w(z)=\sum_{{\bf n}\in C_w}\, 
{z^{n_r}\over n_1\dots n_r}$. 

By Lemma 8, we have:

$$\nsquare(ab^{r-1})=
 a^{r-1}(a+b) +\sum_{1\leq j\leq r-1} a^j(a+b)^{r-j}$$

\noindent Since 
$C_i=C_{a^ic^{r-i}}$ for all $i\leq r-1$, and
$C_r=C_{a^{r-1}c}$, we get

$$L_{\usquare(ab^{r-1})}(z)= 
\sum_{1\leq i\leq r}\,
\sum_{{\bf n}\in C_i} \,\,{z^{-n_r}\over n_1n_2\dots n_r }$$

Therefore,  Lemma 11 follows from Theorem 7. Q.E.D.

\smallskip
Set 

$$C=\{{\bf n}=(n_1\dots n_r)\in{\bf Z}^r\vert 
0<n_1\leq n_2\dots \leq n_r\}.$$ 

\noindent For  ${\bf n}=(n_1\dots n_r)\in C$, define the number 
$a({\bf n})$ as follows. 
If we have $0<n_1<\dots <n_{r-1}$ set $a({\bf n})=k$. Otherwise,
there exist an index $i\leq r-2$ such that
$0<n_1<n_2 \dots n_{i}=n_{i+1}$. In such a case, set
$a({\bf n})=i$.   Note that $a({\bf n})$ does not depend on
the last component $n_r$ of ${\bf n}$, and the function 
${\bf n}\mapsto a({\bf n})$ takes value in the
set $\{1,\,2,\,\dots,\,r-2,\,r\}$

\smallskip
COROLLARY 12:
{\it We have:

$$\zeta(r)= 
\sum_{{\bf n}\in C} a({\bf n}) {2^{-n_r}\over n_1 n_2\dots n_r }$$.}

\smallskip
{\it Proof:}
It is easy to check that $a({\bf n})$ is precisely the number of indices
$i$, $1\leq i\leq r$ such that ${\bf n}$ belongs to $C_i$. Therefore
the formula of Corollary 7 follows from Lemma 11. Q.E.D.

\smallskip
Examples:  For $r=2$, then $a((m,n))=2$ for all
$(m,n)\in C$. Therefore, we get

$$\zeta(2)= 
2 \sum_{0<m\leq n} {2^{-n}\over nm}=2L_2(1/2) +\log^2 2$$

\noindent Accordingly to {\bf [C]}, this formula is due to Euler.

 For $r=5$, we have 
$a((k,l,m,n,p))=1$ if $k=l$,
$a((k,l,m,n,p))=2$ if $k<l= m$  
$a((k,l,m,n,p))=3$ if $k<l< m=n$
and $a((k,l,m,n,p))=5$ if $k<l< m<n$. Therefore, we get
the following expansion for $\zeta(5)$
$$\sum\limits_{0<l\leq m\leq n\leq p} {2^{-p}\over l^2mnp}
+2 \sum\limits_{0<l<m\leq n\leq p} {2^{-p}\over lm^2np}
+3\sum\limits_{0<l<m< n\leq p} {2^{-p}\over lmn^2p}
+5\sum\limits_{0<k<l<m< n\leq p} {2^{-p}\over klmnp}
$$

\vskip1cm
{\bf 4. Geometric interpretation of Theorem 7.}

Theorem 7 provides a combinatorial way to exress any polyzeta value
as the value of a polylogarithmic functions at $1/2$
or  at $\rho^{\pm 1}$. The combinatorics seem 
very intricate: e.g.  the explicit formulas for zeta values
$\zeta(r)$ of Section 3 are  difficult to extend for general polyzeta 
values $\zeta(r_1,\dots,r_k)$.

In this section,  Theorem 7  is reformulated in terms of simple geometric
notions.

\smallskip
(4.1) 
First, the free pro-algebraic group on two
generators $\Gamma$ and its Lie algebra $\hbox{\goth g}$
are defined.

  Let $F$ be the free Lie {\bf Q}-algebra with two generators
$\alpha$ and $\beta$,  let $C^n F$ be its central descending series
and set  $\hbox{\goth g}=\lim\limits_{\leftarrow} F/C^nF$. Since
$F/C^nF$  is a nilpotent Lie algebra, the Campbell-Hausdorf series defines a
structure of algebraic group  on $F/C^nF$, denoted by $\Gamma_n$. 
Then $\Gamma=\lim\limits_{\leftarrow} \Gamma_n$ is a proalgebraic group
(an alternative definition of $\Gamma$ is given in the introduction).
As pro-algebraic varieties,  $\hbox{\goth g}$ and  $\Gamma$ are identical,
and the corresponding isomorphism is denoted by  
$\exp: \hbox{\goth g}\rightarrow \Gamma$.

 Let $F=\oplus_{n\geq 1}\, F_n$ be the grading of $F$ such that $F_1={\bf
Q}\alpha \oplus {\bf Q}\beta$. Then we have 
$\hbox{\goth g}=\prod\limits_{n\geq 1} F_n$, so any 
$x\in\hbox{\goth g}$ can be written as the series
$x=\sum_{i>0} x_i$ where $x_i\in F_i$.  The multiplicative group
${\bf Q}^*$ acts linearly on $\hbox{\goth g}$ as follows:
$t.x=\sum_{i>0}\, t^i\, x_i$, for any $t\in {\bf Q}^*$.

\smallskip
{\it LEMMA 13:} {\it Let $\Phi:\hbox{\goth g}\rightarrow\hbox{\goth g}$
be a morphism of pro-algebraic varieties. Assume that $\Phi$ is
${\bf Q}^*$-invariant and that $\hbox{\rom d} \Phi_0$ is invertible, then
$\Phi$ is an isomorphism.}

\smallskip

{\it Proof:} One can assume that $\hbox{\rom d} \Phi_0$  is the identity.
Then choose a basis  of $F$ consisting of homogenous elements
$(e_n)_{n\geq 1}$ with  $d_n\leq d_m$ if $n<m$, where $d_n$ is
the degree of $e_n$. Accordingly, we have
$\Phi(\sum_{n\geq 1} x_n e_n)= \sum_{n\geq 1} \Phi_n(x) e_n$, where each
$\Phi_n$ is a polynomial in $x=(x_1,\,x_2,\dots)$.  By hypothesis, the
linear part of
$\Phi_n(x)$ is $x_n$  and for any monomial $x_{i_1}\dots x_{i_k}$
 occuring in $\Phi_n(x)$ we have $d_{i_1}+\dots d_{i_k}=d_n$.
It follows that  $\Phi_n(x)-x_n$ depends only on $x_1,\dots,x_{n-1}$, so we
can  write: $\Phi_n(x)=x_n+H_n(x_1,\dots,x_{n-1})$. Since
$\Phi$ is triangular, it is an isomorphism.

\smallskip
(4.2)
There is an isomorphism of Hopf algebras
${\bf Q}[\Gamma]\simeq {\cal H}$, see {\bf [P]}.
A natural group isomorphism
$\psi:\Gamma\rightarrow \hbox{\rom Spec}\,{\cal H}$
is now described.

For  $t\in {\bf Q}$,  define two points $\phi_a(t)$ and 
$\phi_b(t)$ in $\hbox{\rom Spec}\,{\cal H}$ as follows.
Since words $w\in {\cal W}$ are functions on  
$\hbox{\rom Spec}\,{\cal H}$, one needs to evaluate 
$w$ at the points $\phi_a(t)$ and 
$\phi_b(t)$. The rule is as follows:

$w(\phi_a(t))= t^n/n!$ if $w=a^n$ and 
$w(\phi_a(t))=0$ if $b$ occurs in $w$.

$w(\phi_b(t))= t^n/n!$ if $w=b^n$ and 
$w(\phi_b(t))=0$ if $a$ occurs in $w$.

Then it is clear that $\phi_a(t)$ and $\phi_b(t)$ are two one-parameter
groups in $\hbox{\rom Spec}\,{\cal H}$. Since 
$\Gamma$ is freely generated (as a proalgebraic group) by the
two one-parameter groups $\exp\,{\bf Q}\alpha$ and
$\exp\,{\bf Q}\beta$ , the  
isomorphism $\psi$ is prescribed
by the requirements
$\psi(\exp\,t\alpha)=\phi_a(t)$ and  $\psi(\exp\,t\beta)=\phi_b(t)$
for all $t\in {\bf Q}$.

\smallskip
(4.3)  From now on, we identify ${\cal H}$ and ${\bf Q}[\Gamma]$.
Since ${\cal H}={\bf Q}[\Gamma]$, any function
$\phi\in {\bf Q}[\Gamma]$ defines a polylogarithmic function $L_{\phi}(z)$
and the polyzeta value $\zeta(\phi)$ and the numbers
$\zeta^{\pm}(\phi)$. 

\smallskip
(4.4)
The maps $\sigma, \tau:{\cal H}\rightarrow{\cal H}$ are
anti-isomorphisms  of Hopf algebras, and therefore they induce two
anti-isomorphisms  of $\Gamma$ and of its Lie algebra $\hbox{\goth g}$.
These  are again denoted by $\sigma$ and $\tau$.
They are  uniquely characterized by the requirements:

\centerline{$\sigma\exp\,t\alpha=\exp\,t\beta$ and 
$\sigma\exp\,t\beta=\exp\,t\alpha$}

\centerline{$\tau\exp\,t\alpha=\exp\,t(\alpha+\beta)$ and 
$\tau\exp\,t\beta=\exp\,-t\alpha$,}

for all $t\in{\bf Q}$. We have $\sigma^2(g)=g$ and 
$\tau^3(g)=g^{-1}$ for any $g\in G$.

\smallskip
(4.5) Since ${\cal H}={\bf Q}[\Gamma]$,
the maps $\square$, $\nabla$ occuring in Theorem 7  are now identified with
some algebra morphisms
$\square, \nabla: {\bf Q}[\Gamma]\rightarrow{\bf Q}[\Gamma]$.

\smallskip
LEMMA 14: {\it Let $\phi\in {\bf Q}[\Gamma]$. Then for any
$g\in\Gamma$, we have

$$\nsquare\,\phi(g)=\phi(g\sigma(g))\hskip.5cm
\hbox{\rom and}\hskip.5cm
\nabla\,\phi(g)=\phi(g\tau(g))$$}

\smallskip
{\it Proof:}
 Using their definitions, $\rsquare$ and $\nabla$ are the composition
of the following maps:

$$\Gamma\buildrel\hbox{\srom diag}\over\rightarrow
\Gamma\times\Gamma\buildrel id\times\sigma\over\longrightarrow
\Gamma\times\Gamma\buildrel \mu\over\rightarrow
\Gamma$$

$$\Gamma\buildrel\hbox{\srom diag}\over\rightarrow
\Gamma\times\Gamma\buildrel id\times\tau\over\longrightarrow
\Gamma\times\Gamma\buildrel \mu\over\rightarrow
\Gamma$$

where $\hbox{\rom diag}(g)=(g,g)$ and $\mu(g_1,g_2)=g_1.g_2$.
Therefore we have 

\centerline{$\nsquare\,\phi(g)=\phi(g.\sigma(g))$
and
$\nabla\,\phi(g)=\phi(g.\tau(g))$.}

\smallskip
(4.6) In this subsection, the symmetric space 
associated with $\sigma$ is defined.

Set 
$\hbox{\goth k}=\{x\in \hbox{\goth g}\vert \sigma(x)=-x\}$,
$K=\{g\in \Gamma\vert \sigma(g)=g^{-1}\}$, 
$\hbox{\goth p}=\{x\in \hbox{\goth g}\vert \sigma(x)=x\}$
$P=\{p\in \Gamma\vert \sigma(g)=g\}$.

Since $\sigma$ is an anti-involution,  $K$ is a subgroup in $\Gamma$
and $\hbox{\goth k}$ is its a Lie algebra. Obviously we have
$\hbox{\goth g}=\hbox{\goth k}\oplus\hbox{\goth p}$.
Since $\Gamma$ is a pro-unipotent group, we have
$K=\exp\,\hbox{\goth k}$, $P=\exp\,\hbox{\goth p}$
and $\Gamma=P.K$. So any element $g\in \Gamma$ can be written
as $g=p.k$, where $k\in K$ and $p\in P$. Moreover $P\simeq \Gamma/K$ is a
symmetric space.

\bigskip
LEMMA 15: {\it Let $\phi\in {\bf Q}[\Gamma]$. Then for any
$g=p.k\in\Gamma$, we have

$$\nsquare\,\phi(p.k)=\phi(p^2)$$

In particular, $\zeta(\phi)=0$ if $\phi\vert_P\equiv 0$. }

\smallskip
{\it Proof:} This follows from Lemma 14 and Theorem 7.

\smallskip
(4.7) Note  that $\tau$ is not an involution, but an "ordrer three"
anti-isomorphism, i.e. $\tau^3(g)=g^{-1}$. In this sub-section,
we introduce a
space $Q$ which is analogous to a symmetric space.

Set $\hbox{\goth l}=\{x\in \hbox{\goth g}\vert \tau(x)=-x\}$,
$L=\{g\in \Gamma\vert \tau(g)=g^{-1}\}$, 
$\hbox{\goth q}=\{x\in \hbox{\goth g}\vert 
\tau^2(x)-\tau(x)+x=0\}$. Also define
$Q$ as the image of the map $g\in\Gamma\mapsto g\tau(g)$.
Note that $L$ is a subgroup with Lie algebra 
$\hbox{\goth l}$.

\smallskip
{\it LEMMA 16:} {\it The subset $Q$ is a closed subvariety of $\Gamma$
and the natural map: $Q\times L\rightarrow\Gamma,\,(q,l)\mapsto ql$
is a isomorphism of pro-algebraic varieties.}

\smallskip
{\it Proof} It is easy to prove
that $\Gamma=\exp\hbox{\goth q}.L$. For 
$g= \exp q.l$, with $q\in \hbox{\goth q}$ and $l\in L$, we have
$g\tau(g)=\exp q \exp\tau(q)$, therefore $Q$ is the set all
$\exp q \exp\tau(q)$ for $q\in \hbox{\goth q}$.

Let $\Phi:\hbox{\goth q}\oplus \hbox{\goth l}\rightarrow \Gamma$
be define by $\Phi(q,l)=\exp q \exp\tau(q)\exp l$. 
Note that $\hbox{\rom d} \Phi_0$ is the linear map
from $\hbox{\goth g}$ to $\hbox{\goth g}$ which is the identity on
$\hbox{\goth l}$ and whose restriction to
$\hbox{\goth q}$ is $1+\tau$. Therefore, $\hbox{\rom d} \Phi_0$
is invertible. By Lemma 13 that $\Phi$ is an isomorphism
and  Lemma 16 follows easily.
Q.E.D.

\smallskip
The definition of $Q$ is slighty more complicated than the definition
of $P$ because $Q\neq \exp\hbox{\goth q}$. However, the
map $\hbox{\goth q}\rightarrow Q,\, q\mapsto \exp q\,\exp\tau(q)$
is an isomorphism from $\hbox{\goth q}$ to $Q$.
Any element $g\in \Gamma$ can be written
as $g=q.l$, where $l\in L$ and $q\in Q$

\bigskip
LEMMA 17: {\it Let $\psi\in {\bf Q}[\Gamma]$. Then for any
$g=q.l\in\Gamma$, we have

$$\nabla\,\psi(q.l)=\psi(q\tau(q))$$

In particular, $\zeta(\psi)=0$ if $\psi\vert_Q\equiv 0$.}

\smallskip
{\it Proof:} This follows from Lemma 14 and Theorem 7.

\smallskip
(4.8) Let $\phi\in{\bf Q}[P]$ be a rational function on $P$.
The notations $\zeta(\phi)$ and $L_{\phi}(z)$ are now defined.
Set  $\zeta(\phi)=\zeta(\hat\phi)$ where
$\hat\phi$ is any function on $\Gamma$ extending $\phi$. By Lemma
15, $\zeta(\phi)$ is well defined.
Since  $P\simeq \Gamma/K$,
$\phi$ can be uniquely extended to  a right $K$-invariant function $\Phi$
on $\Gamma$. Then set  $L_{\phi}(z)=L_{\Phi}(z)$. 

Similarly, for $\psi\in{\bf Q}[Q]$,
the notations $\zeta^{\pm}(\psi)$ and $L_{\psi}(z)$ are defined as follows.
Set  $\zeta^{\pm}(\psi)=\zeta^{\pm}(\hat\psi)$ where
$\hat\psi$ is any function on $\Gamma$ extending
$\psi$. By Lemma
17, $\zeta^{\pm}(\psi)$ is well defined. 
By Lemma 16, we have  $Q\simeq \Gamma/L$, therfore $\psi$
can be uniquely extended to  a right $L$-invariant function $\Psi$ on
$\Gamma$. Then set 
$L_{\psi}(z)=L_{\Psi}(z)$.

Define the algebra morphisms
$\rsquare:{\bf Q}[P]\rightarrow {\bf Q}[P]$
and $\nabla:{\bf Q}[Q]\rightarrow {\bf Q}[Q]$ by

$$\nsquare\phi(p)=\phi(p^2), 
\hskip0.4cm\hbox{\rom for} \phi\in {\bf Q}[P]$$

$$\nabla\psi(q)=\psi(q\tau(q)), 
\hskip0.4cm\hbox{\rom for} \psi\in {\bf Q}[Q]$$

These operators are simply
the restrictions to $P$ and to $Q$ of the
already defined operators 
$\rsquare, \nabla: {\bf Q}[\Gamma]\rightarrow {\bf Q}[\Gamma]$.  So
using the same notations should  not bring  confusions.

\bigskip
THEOREM 18: {\it For any $\phi\in {\bf Q}[P]$ and
$\psi\in {\bf Q}[Q]$, we have

$$\zeta(\psi)= L_{\msquare(\psi)}(1/2)\hskip0.4cm\hbox{\rom
and}\hskip0.4cm\zeta^{\pm}(\psi)= L_{\nabla(\psi)}(\rho^{\pm})$$}

\bigskip {\it Proof:} It follows immediately from Theorem 7, and
Lemmas 15 and 17.

\vskip1cm

{\bf 5. Other expressions for zeta values. }

In this section, we follow a suggestion of W. Zudilin

\smallskip
(5.1) Theorem 7 shows that any polyzeta value is the value of a 
polylogarithmic function at $1/2$ or at $\rho^{\pm}$. However, there is a
much more simple way to express the zeta values $\zeta(r)$ as a   value of
polylogarithmic functions at
$1/2$ or at $\rho^{\pm}$, see Corollary 20.
 It is surprizing that the two approaches give different expressions, except
for $\zeta(2)$. Moreover, this simpler approach does not generalize to
polyzeta values.

\smallskip
(5.2)
Let $\sigma':{\cal H}\rightarrow {\cal H}$ be the 
linear map defined as follows. Set $\sigma'(a)=-a$ and
$\sigma'(b)= a+b$. For a word $w=c_1\dots c_n$, where 
$c_i\in\{a,\,b\}$ are letters, set 
$\sigma'(w)=\sigma'(c_1)\dots \sigma'(c_n)$. It is easy to prove 
that $\sigma'$ is an algebra morphism relative to the schuffle product
and that $\sigma'({\cal H^+})={\cal H^+}$.

\smallskip
LEMMA 19 {\it  For any $h\in {\cal H^+}$, we have:

\centerline{$L_h(-1)=L_{\sigma'(h)}(1/2)$ and
$L_h(\overline\rho)=L_{\sigma'(h)}(\rho)$. }}

\smallskip
{\it Proof:} One can assume that $h$ is a word
$w=c_1\dots c_n\in{\cal W}^+$. Set 
$F=\{z\in{\bf C}\vert\,\vert z\vert\leq 1
\,\hbox{\rom and}\, Im z\leq 1/2\}$ and 
choose a two paths $\gamma,\,\gamma_-:[0,1]\rightarrow F$
with $\gamma(0)=\gamma_-(0)=0$, $\gamma(1)=-1$
and $\gamma_-(1)=\overline\rho$.
By Proposition 2, we have:

 $$L_{w}(-1)=\int_{\Delta_n}\,\gamma^*\omega_{c_1}(x_1)
\gamma^*\omega_{c_2}(x_2)\dots\gamma^*\omega_{c_n}(x_n)$$

$$L_{w}(\overline\rho)=\int_{\Delta_n}\,\gamma_-^*\omega_{c_1}(x_1)
\gamma_-^*\omega_{c_2}(x_2)\dots\gamma_-^*\omega_{c_n}(x_n)$$

For $z\in F$, set $S'(z)=z/(z-1)$ and  
set $\delta=S'\circ\gamma$ and $\delta_-=S'\circ\gamma_-$.
We have $S'^*\omega_{\sigma'(c)}=\omega_c$ for $c=a$ or $b$, 
and therefore we get
$\gamma^*\omega_c=\delta_-^*\omega_{\sigma'(c)}$ and 
$\gamma_-^*\omega_c=\delta_-^*\omega_{\sigma'(c)}$ 
It follows that

 $$L_{w}(-1)=\int_{\Delta_n}\,\delta^*\omega_{\sigma'(c_1)}(x_1)
\delta^*\omega_{\sigma'(c_2)}(x_2)\dots\delta^*\omega_{\sigma'(c_n)}(x_n)$$

$$L_{w}(\overline\rho)=\int_{\Delta_n}\,\delta_-^*\omega_{\sigma'(c_1)}(x_1)
\delta_-^*\omega_{\sigma'(c_2)}(x_2)\dots\delta_-^*\omega_{\sigma'(c_n)}(x_n)$$

We have  $S'(F)=F$, $S'(0)=0$,
$S'(-1)=1/2$ and $S'(\overline\rho)=\rho$. Therefore 
$\delta$ is a path from $0$ to $1/2$ and
$\delta_-$ is a path from $o$ to $\rho$.
Thus, these integrals
can be identified  by Proposition 2, and we get
$L_h(-1)=L_{\sigma'(h)}(1/2)$, and
$L_h(\overline\rho)=L_{\sigma'(h)}(\rho)$. Q.E.D.

\smallskip
{\it LEMMA 20}: {\it  Let $r\geq 1$.  We have

$$\zeta(r+1)={-1\over1-2^{-r}}\,L_{r+1}(-1)$$

$$\zeta(r+1)={1\over(1-2^{-r})(1-3^{-r})}\,\,
[L_{r+1}(\rho)+L_{r+1}(\overline\rho)]$$}

\smallskip
{\it Proof:} For each positive integer $a$, set
$\delta_a(n)=1$ if $a$ divides  $n$  and 
$\delta_a(n)=0$ otherwise.

From the formula $(-1)^n=-\delta_1(n)+2\delta_2(n)$, we get

$L_{r+1}(-1)=\sum\limits_{n>0} (-1)^n/n^{r+1}$

\hskip1.65cm$=-\sum\limits_{n>0} 1/n^{r+1} + 2\sum\limits_{n>0} 1/(2n)^{r+1}$

\hskip1.65cm$=(-1+2^{-r})\,\zeta(r)$

\noindent from which the first formula follows.

From the formula $\rho^n+\overline{\rho}^n=
\delta_1(n)-2\delta_2(n)-3\delta_3(n)+6\delta_6(n)$
we get 

$L_{r+1}(\rho)+L_{r+1}(\overline\rho)=
\sum\limits_{n>0} [\rho^n +\overline\rho^n]/n^{r+1}$

\hskip2.95cm$=\sum\limits_{n>0} 1/n^{r+1}- 2\sum\limits_{n>0} 1/(2n)^{r+1}
-3\sum\limits_{n>0} 1/(3n)^{r+1}+6\sum\limits_{n>0} 1/(6n)^{r+1}$

\hskip2.95cm$=(1-2^{-r}-3^{-r}+6^{-r})\,\zeta(r+1)$

\hskip2.95cm$=(1-2^{-r})(1-3^{-r})\,\zeta(r+1)$,

\noindent  from which the second formula follows.

\smallskip
COROLLARY 21:
{\it Let $r\geq 1$.  We have

$$\zeta(r+1)= {2^r\over (2^r-1)}\,L_{a(a+b)^r}(1/2)$$

$$\zeta(r+1)= {6^r\over (2^r-1)(3^r-1)}\,
[L_{ab^r}(\rho)-L_{a(a+b)^r}(\rho)]$$
}

\smallskip
{\it Proof:} The corollary follows from Lemmas 18 and 19.

\bigskip
Examples:

$$\zeta(2)= 
2 \sum_{0<m\leq n} {2^{-n}\over nm}$$

$$\zeta(3)= 
8/7 \sum_{0<l\leq m\leq n} {2^{-n}\over lnm}$$

The expression for $\zeta(2)$ is the
same as in section 3. However for all other zeta values $\zeta(r)$ with
$r\geq 3$, the expressions are different: e.g., the formula of
corollary 20 uses non integral coefficients. Moreover, this simpler
approach only concerns zeta values but not the polyzeta values.

\vskip1cm

{\bf 6. Conclusion:} Polyzeta values are mixed periods
{\bf [G]}, {\bf [T]}, {\bf [Z]}.
In the philosophy of motives, there is a proalgebraic 
group $G$ and  periods should be regular functions on $G$:
more precisely, the algebra of periods should be
a $\overline{\bf Q}$-form of ${\bf C}[G]$, modulo conjectures. 
Here polyzeta
values are attached to  some symmetric space. Does there is a
motivic interpretation of this construction? Note however that
the map ${\bf Q}[P]\rightarrow {\cal Z}$ is not injective.

\vskip1cm

{\it Bibliography:}

\smallskip
{\bf [BBP]}
Bailey, David,  Borwein, Peter , and Plouffe, Simon: On
the Rapid Computation of Various Polylogarithmic Constants.
Mathematics of Computation 66 (1997)  903-913.

\smallskip
{\bf [C]}
P. Cartier: Polylogarithmes, polyzetas et groupes pro-unipotents.
{\it S\'eminaire Bourbaki 2000-2001} Ast\'erisque 282 (2002) 137-173

\smallskip
{\bf [E]}
Ecalle, Jean: La libre g\'en\'eration des multiz\^etas et leur
d\'ecomposition canonici-explicite en irr\'eductibles, Preprint
(1999).

\smallskip
{\bf [G]}  Goncharov, Alexander:
Multiple $\zeta$-values,
Galois groups, and geometry of modular varieties. 
{\it European Congress of
Mathematics (Barcelona, 2000).} Progr. Math. 201 (2001)
361-392. 

\smallskip
{\bf [H]}
M.E. Hoffman: Multiple harmonic series, J. of Algebra 194
(1992) 275-290

\smallskip
{\bf [L]}
Lewin, Leonard: {\it Polylogarithms and Associated Functions.}
North-Holland (1981)

\smallskip
{\bf [MP]}
Minh, Hoang Ngoc,  and Petitot, Michel Lyndon: Words, 
polylogarithms and the Riemann $\zeta$ function. Discrete Math. 217
(2000) 273--292.

\smallskip
{\bf [O]}
Oesterl\'e, Joseph:  Polylogarithmes. {\it S\'eminaire
Bourbaki}.  Ast\'erisque  216  (1993) 49-67.

\smallskip
{\bf [P]}
Pianzola, Arturo:  Free group functors.  J. Pure Appl. Algebra  140 
(1999) 289-297.

\smallskip
{\bf [T]} Terasoma, Tomohide: Mixed Tate motives and multiple zeta values.
Invent. Math. 149 (2002) 339-369.

\smallskip
{\bf [Z]}
Zagier, Don:  Values of zeta functions and their applications. First European
Congress of Mathematics (Paris, 1992), 
Birkh\"auser, Basel, Progr. Math., 120 (1994) 497-512.

\vskip1cm
{\it Author address:}

mathieu@math.univ-lyon1.fr

Universit\'e de Lyon

Institut Camille Jordan, UMR 5028 du CNRS

43 bd du 11 novembre 1918

69622 Villeurbanne cedex 

France

\end